\documentclass[final,epsfig,amsfonts]{article}
\usepackage{latexsym}
\usepackage[dvips]{graphicx}
\newcommand{\eeq}{\end{equation}}
\newcommand{\beq}{\begin{equation}}
\newcommand{\nuq}[1]{\label{#1} \eeq}
\newcommand{\ovl}[1]{\overline{#1}}

\newtheorem{definition}{Definition}

\newtheorem{proposition}{Proposition}
\newtheorem{remark}{Remark}
\newtheorem{lemma}{Lemma}

\newtheorem{criterion}{Criterion}

\newtheorem{theorem}{Theorem}

\def\bfq{\mathbf{Q}}
\def\bfn{\mathbf{N}}
\def\bfz{\mathbf{Z}}

\begin{document}
\title{
Minkowski's question mark measure is UST--regular}
\author{
Giorgio Mantica \\
Center for Non-linear and Complex Systems, \\ Dipartimento di Scienza ed Alta Tecnologia, \\Universit\`a dell'~Insubria, 22100 Como, Italy \\ and \\ CNISM unit\`a di Como,
INFN sezione di Milano, \\
Istituto Nazionale di Alta Matematica, \\
Gruppo Nazionale per la Fisica Matematica.}

\date{}
\maketitle
\begin{abstract}
We prove the recent conjecture that Minkowski's question mark measure is regular, in the sense of Ullman--Stahl--Totik.

{\em keywords: Minkowski question mark function; Regular Measures; M\"obius Iterated Function Systems }\\
{\em MATH Subj. Class. 42C05; 31A15; 11B57; 28A80; 37A45}

\end{abstract}

\section{Introduction and statement of the main results}

A remarkable function has been introduced by Hermann Minkowski in 1904, to map  algebraic numbers of second degree to the rationals (and these latter to binary fractions) in a continuous, order preserving way \cite{minko}. This function is called the {\em question mark function} and it is indicated by $?(x)$, perhaps because of its enigmatic yet captivating, multi--faceted personality. In fact, it is linked to continued fractions, the Farey tree and the theory of numbers. It also appears in the theory of dynamical systems, in relation with the Farey shift map \cite{lagarias,isola,mirko} and in the coding of motions on manifolds of negative curvature \cite{series,guzzi,guzzi-ke,prlmob,kesse2}.

Let us therefore briefly introduce this function.
Consider the interval $I=[0,1]$ and let $x \in I$. Write this latter in its continued fraction representation, $x = [n_1,n_2,\ldots,]$, set
$N_j(x) = \sum_{l=1}^j n_l$, and define $?(x)$ as the sum of the series \cite{denjoy,salem}
\beq
 ?(x) = \sum_{j=1}^\infty (-1)^{j+1} 2^{-N_j(x)+1}.
\nuq{eq-minko1}
To deal with rational values $x \in I$, also stipulate that terminating continued fractions correspond to finite sums in the above series.

The analytical properties of the question mark function are so interesting that its graph has been named the {\em slippery devil's staircase} \cite{guzzi}: it is continuous and H\"older continuous of order $\log 2/(1+\sqrt{5})$ \cite{salem}. It can be differentiated almost everywhere; its derivative is almost everywhere null \cite{denjoy,salem} and yet it is strictly increasing: $?(y)-?(x) > 0$ for any $x,y \in I$, $x<y$. The fractal properties of the level sets of its derivative have been studied via the multifractal formalism \cite{guzzi,kesse2}.

Since $?(x)$ is monotone non--decreasing, it is the distribution function of a Stieltjes measure $\mu$:
\begin{equation}
  ?(x) = \mu([0,x)),
\label{eq:defmu}
\end{equation}
which, because of the above, turns out to be singular continuous with respect to Lebesgue. We call $\mu$ the Minkowski's question mark measure.
A remarkable result by Kinney \cite{kinney} asserts that its Hausdorff dimension can be expressed in terms of the integral of the function $\log_2(1+x)$ with respect to the measure $\mu$ itself. Very precise numerical estimates of this dimension have been obtained with high precision arithmetics \cite{alkatesi}; rigorous numerical lower and upper bounds \cite{minkio} derived from the Jacobi matrix of $\mu$ place this value between 0.874716305108213 and 0.874716305108207. Further analytical properties of $\mu$ have been recently studied, among others, in \cite{alau1,alau2,paradis}.

In this paper we are concerned with different, yet related, fine properties of Minkowski's question mark measure, stemming from logarithmic potential theory in the complex plane \cite{ran0,saff}. In this context, Dresse and Van Assche \cite{dresse} asked whether it is regular, in the sense of Ullmann--Stahl--Totik, as defined below. Their numerical investigation was successively refined via a more powerful technique by the present author in \cite{minkio}, to provide compelling numerical evidence in favor of regularity of this measure. In this paper we provide a rigorous proof of this result, which further reveals the intriguing nature of Minkowski's question mark function.

The notion of {\em regularity} of a measure, due to Ullman \cite{ullman}, Stahl and Totik \cite{stahl}, concerns the asymptotic properties of its orthonormal polynomials $p_j(\mu;x)$ (recall the defining property: $\int p_j(\mu;x) p_m(\mu;x) d \mu(x) = \delta_{jm}$, where $\delta_{jm}$ is the Kronecker delta). In the present case, where the support of the measure $\mu$ is the full interval $[0,1]$, which has capacity $1/4$ and carries a unique equilibrium measure, regularity of $\mu$ means that its orthogonal polynomials somehow mimic Chebyshev polynomials---that are orthogonal with respect to the equilibrium measure---both in {\em root asymptotics} away from $[0,1]$ and in the asymptotic distribution of their zeros in $[0,1]$. Formally, letting $\gamma_j$ be the (positive) coefficient of the highest order term, $p_j(\mu;x) = \gamma_j x^j + \ldots$, regularity is defined by the fact that the geometric mean $\gamma_j^{1/j}$ tends to the capacity $1/4$, when the order $j$ tends to infinity. Equivalent definitions of regularity can be found in \cite{stahl}, collected in definition 3.1.2. A wealth of potential--theoretic properties follow from regularity, as discussed in detail in Chapter 3 of \cite{stahl}.

In the case of Minkowski's question mark measure, regularity and the asymptotic behavior of orthogonal polynomials have been investigated theoretically and numerically in \cite{minkio}, with detailed pictures illustrating the abstract properties. This investigation continues in this paper from a slightly different perspective: we do not prove regularity of $\mu$ directly from the definition, that is, orthogonal polynomials play no r\^ole herein, but we use a purely measure--theoretical criterion, which translates the idea that a regular measure is never too thin on its support. This is Criterion $\lambda^*$ of Stahl and Totik:  \cite{stahl}, Thm. 4.2.7. It reads as follows: if the support of $\mu$ is $[0,1]$ and if for every $\eta>0$ the Lebesgue measure of
\beq
  \{ x \in [0,1] \mbox{ s.t. } \mu([x-1/n,x+1/n]) \geq e^{-\eta n} \}
\nuq{eq-def1}
tends to one, when $n$ tends to infinity, then $\mu$ is regular. Our fundamental result is therefore
\begin{theorem}
Minkowski's question mark measure satisfies criterion $\lambda^*$ and hence is regular.\label{teo-1}
\end{theorem}
Let us now describe the steps of the proof of this result and then place it into wider perspective.

Firstly, we need a more transparent definition of Minkowski's question mark function, than offered by eq. (\ref{eq-minko1}): this is provided by the symmetries of $?(x)$, which permit to regard it as the invariant of an Iterated Function System (IFS) \cite{hut,dem} composed of M\"obius maps, following \cite{prlmob,mobius}. We review this approach in Section \ref{sec-1}. In {\bf Lemma \ref{lem-1a}} we show how it can be used to {\em define a countable family of partitions of $[0,1]$ in a finite number of intervals, with elements labeled by words in a binary alphabet}. The notable characteristic of any of these partitions is that all its elements have the same $\mu$--measure---and obviously different length. We develop all remaining theory from this approach, proving every detail of the procedure, so that the paper is fully self--contained and the reader has no need of external material. Of course, a large body of knowledge already exists on these topics. This is true in particular for the relation of Minkowski's question mark function to the Farey tree and Stern--Brocot sequences. In Section \ref{sec-2}, {\bf Lemma \ref{lem-3}}, we show that {\em these sequences coincide with the ordered set of end--points in the M\"obius IFS partitions of $[0,1]$}. As mentioned, none of these results is totally new, but we present them in a concise set--up, that of IFS, which is both elegant and renders sequent analysis easier.

In Section \ref{sec-3} the above theory permits to {\em single out a family of intervals within an IFS partition whose points fulfill the $\lambda^*$ criterion}, when this latter is properly rescaled to fit in the IFS set--up: this is {\bf Criterion \ref{cri-one}} and the result is {\bf Lemma \ref{lem-4}}. Our goal is then to prove that {\em the (Lebesgue) measure of such family tends to one, when the index of the partition tends to infinity}. 
The proof of this fact requires a fine control of the intervals composing M\"obius IFS partitions---equivalently, Stern--Brocot intervals. In fact, we define the set of ``large'' intervals as those whose length is larger than $\alpha/n$, where $\alpha$ is an arbitrarily small number and $n$ is the order of the IFS partition to which they belong. 
In {\bf Proposition \ref{prop-1}} we prove that, {\em for any real positive $\alpha$, the cardinality of large intervals is superiorly bounded, independently of $n$}. 
This interesting result is loosely related to the pressure function from the so--called thermodynamical formalism, that gauges the exponential growth rate of sums of the partition interval lengths, raised to a real power. These sums, for Stern--Brocot intervals, have been studied in \cite{roberto,kesse}.
The paper is then concluded by the proof of {\bf Theorem \ref{teo-1}}, which at this point is slightly more than a judicious exploitation of the previous results.

The fact that Minkowski's question mark function is UST-regular is remarkable in many ways. First, it was not at all obvious how to reveal it from the numerical point of view: it required dedicated techniques \cite{minkio}. From the theoretical side, regularity 
is required in the hypotheses of Proposition 1 and 2 of \cite{minkio}, that are therefore now rigorously established: these propositions describe and quantify in a precise way the {\em local} asymptotic behavior of zeros of the orthogonal polynomials $p_n(\mu;x)$ and of the Christoffel functions associated with $\mu$, and they link them to the Farey / Stern--Brocot organization of the set of rational numbers.

Further conjectures were presented in \cite{minkio}, on the speed of convergence in the above asymptotic behaviors and, more significantly, on the fact that Minkowski's question mark might belong to Nevai's class: numerical indication is that its Jacobi matrix elements converge to a limit value, although slowly. If confirmed, this conjecture will provide us with an example of a measure in Nevai's class which does not fulfill Rakhmanov sufficient condition \cite{rakh,totmany} (almost everywhere positivity of the Radon Nikodyn derivative of $\mu$ with respect to Lebesgue) and might perhaps indicate a widening of such condition; it is known that Nevai's class does contain pure point \cite{alphwalt} and singular measures \cite{doron} but these examples do not seem to indicate a general criterion on a par with Rakhmanov's.

In conclusion, the picture of Minkowski's question mark measure that emerges from recent investigations is that of a singular continuous measure that nonetheless has many {\em regular} characteristics: it is regular according to Ullman--Stahl--Totik; we conjectured that it belongs to Nevai's class; its Fourier transform tends to zero polynomially \cite{yakubo1,yakubo2,persson,jordan} even if it does not fulfill the Riemann--Lebesque sufficient condition. It is therefore interesting---and a direction of further research---to study the so--called Fourier--Bessel functions \cite{etna} generated by Minkowski's question mark measure, in order to detect whether they display any of those features associated with singular continuous measures \cite{str1,str2,str3,igm,gioitalo,guzze} that are typical of measures with almost--periodic Jacobi matrices \cite{physd1,physd2,etna}---which the present measure is conjectured not to be.

\section{Minkowski's question mark measure and M\"obius IFS}
\label{sec-1}
In our view, the most effective representation of Minkowski's question mark function is via an Iterated Function System \cite{dem,hut} composed of M\"obius maps \cite{prlmob,mobius}. This is a translation in modern language of the relation between Minkowski's question mark function and modular transformations already discussed in \cite{denjoy}. Let us therefore adopt and develop the formalism introduced in \cite{prlmob}.
Define maps $M_i$ and $P_i$, $i=0,1$ from $[0,1]$ to itself as follows:
\begin{equation}
\label{eq-mink1}
    \begin{array}{ll}
    M_0(x) = \frac{x}{1+x}, & P_0(x) = \frac{x}{2},\\
    M_1(x) = \frac{1}{2-x}, & P_1(x) = \frac{x+1}{2}
    \end{array}
\end{equation}
Then, using the properties of the continued fraction representation of a real number and eq. (\ref{eq-minko1}) (see {\em e.g.} \cite{prlmob}) it is not difficult to show that the following properties hold :
\begin{equation}
\label{eq-mink1a}
    ?(0) =  0, \;\; ?(1) = 1,
\end{equation}
\begin{equation}
\label{eq-mink1b}
    ?(M_i(x)) =  P_i(?(x)), \; i = 0,1.
\end{equation}
We proved in \cite{prlmob,mobius} that these relations uniquely define the function $?(x)$, since they define an Iterated Function Systems, composed of the two M\"obius maps $M_i$, whose invariant measure is Minkowski's question mark measure $\mu$. We will not need the full power of this construction, which has been exploited also in \cite{minkio}; nonetheless we will use equations (\ref{eq-mink1})--(\ref{eq-mink1b}) as a basis of all our theory. We start from the following
\begin{definition}[Symbolic Words]
Let $\Sigma$ be the set of finite words in the letters $0$ and $1$.  We denote by $|\sigma|$ the {length} of $\sigma \in \Sigma$: if  $|\sigma|=n$ then
$\sigma$ is the $n$-letters sequence $\sigma_1,\sigma_2,\ldots,\sigma_n$ where $\sigma_i$ is either $0$ or $1$. Let $\emptyset$ be the empty word and assign to it length zero. Denote by $\Sigma^n$ the set of $n$-letter words, for any $n \in \bfn$. Given two words $\sigma \in \Sigma^n$ and $\eta \in \Sigma^m$ the composite word $\sigma \eta \in \Sigma^{n+m}$ is the sequence $\sigma_1,\ldots,\sigma_n,\eta_1,\ldots,\eta_m$.
\label{def-1}
\end{definition}
\begin{lemma}
\label{lem-1a}
Let $\Sigma^n$ be as in definition \ref{def-1}. Associate to any $\sigma \in \Sigma^n$ the map composition
\[
 M_\sigma = M_{\sigma_1} \circ M_{\sigma_2} \circ \cdots \circ M_{\sigma_n},
\]
when $n >0$, and let $M_\emptyset$ be the identity transformation. Let
\[I_\sigma= M_\sigma([0,1]).\]
Then, for any integer value $n \in \bfn$, the intervals $I_\sigma$, with $\sigma \in \Sigma^n$, are pairwise disjoint except possibly at one  endpoint and fully cover $I$:
\beq
 I = \bigcup_{\sigma \in \Sigma^n} I_{\sigma}.
\nuq{eq-ifull}
 \end{lemma}
{\em Proof.}
When $\sigma=\emptyset$ the lemma is obvious.
Observe that the functions $M_i$, $i=0,1$ are continuous, strictly increasing and map $[0,1]$ to the two intervals $[0,\frac{1}{2}]$ and $[\frac{1}{2},1]$ respectively, which are disjoint except for a common endpoint. Then, the same happens for the two intervals $(M_{\sigma} \circ M_i) ([0,1]) = I_{\sigma i}$, $i=0,1$, where $\sigma$ is any finite word and  $\sigma i$ is the composite word. Explicit computation yields
\[
 I_{\sigma 0}  = [ M_\sigma (M_0 (0)), M_\sigma (M_0(1))]
 = [ M_\sigma (0), M_\sigma (\frac{1}{2})]
\]
and
\[
 I_{\sigma 1}  = [ M_\sigma (M_1 (0)), M_\sigma (M_1(1))]
 = [ M_\sigma (\frac{1}{2}), M_\sigma (1)],
\]
where we have used a property that will be useful also in the sequel: for any $\sigma \in \Sigma$
\beq
M_{\sigma 0} (1) = M_{\sigma 1} (0) = M_\sigma(\frac{1}{2}),
\nuq{eq-mchange}
which is valid since $M_1(0)=M_0(1)=1/2$.
It follows from this that $I_{\sigma 0}$ and $I_{\sigma 1}$ not only are adjacent, but also they exactly cover $I_\sigma$:
\beq
 I_{\sigma 0} \bigcup I_{\sigma 1} = I_\sigma.
\nuq{eq-cover}
Using induction one then proves eq. (\ref{eq-ifull}). $\Box$ \\
As a consequence of this Lemma, each set $\Sigma^n$ is associated with a partition of $[0,1]$ produced by the M\"obius IFS. Since any word in $\Sigma^n$ is uniquely associated to an interval of this partition, in the text we will use the terms word and interval as synonyms.

\begin{lemma}
\label{lem-2a}
Let $\Sigma^n$ be as in definition \ref{def-1}. For any $n \in \bfn$ the function
\beq
\Theta(\sigma) = \sum_{j=1}^{n} \sigma_j 2^{n-j},
\nuq{eq-order1}
induces the lexicographical order in $\Sigma^n$, in which the letter 1 follows the letter 0 and we read words from left to right. In addition, letting
\beq
x_\sigma = M_\sigma(0) = M_{\sigma_1} \circ  \cdots \circ  M_{\sigma_n} (0)
\nuq{eq-order10}
the set $\{x_\sigma, \sigma \in \Sigma^n \}$ is increasingly ordered: $x_\sigma < x_\eta$ if and only if $\sigma<\eta$.
Finally, one has that
\beq
I_\sigma = [x_\sigma,x_{\hat \sigma}]
\nuq{eq-order1b}
where $\hat \sigma$ is the successive word of $\sigma$ when $\sigma \neq 1^n$ and $x_{\hat \sigma} = 1$ in the opposite case.
\end{lemma}
{\em Proof.}
Observe that when $n=0$ we have $\sigma=\emptyset$ and $\Theta(\sigma)=0$ because the sum in (\ref{eq-order1}) contains no terms.
It is immediate that $\Theta$ is bijective from $\Sigma^n$ to $\{0,\ldots,2^n-1\}$ and therefore it induces an order on $\Sigma^n$. This coincides with the lexicographical order that we denote by '$<$': in fact $\sigma < \eta$ if and only if $\Theta(\sigma) < \Theta(\eta)$. To prove this statement, if $\sigma \neq \eta$ we can define $k = \min \{ j \mbox{ s.t. } \sigma_j \neq \eta_j\}$. Then, $\sigma < \eta$ happens if and only if $\sigma_k=0$ and $\eta_k=1$. But in this case one has that
\[
\Theta(\sigma) = \sum_{j=1}^{k-1} \sigma_j 2^{n-j} + 0  + \sum_{j=k+1}^{n} \sigma_j 2^{n-j}
\]
and
\[
\Theta(\eta) = \sum_{j=1}^{k-1} \eta_j 2^{n-j} +  2^{n-k} + \sum_{j=k+1}^{n} \eta_j 2^{n-j}.
\]
In the above equations, the first summations at r.h.s. are equal, since $\sigma_j=\eta_j$ for $j<k$. In addition, the second summation in $\Theta(\sigma)$ is strictly less than $2^{n-k}$ for any choice of the sequence $\sigma_{k+1},\ldots,\sigma_n$ and therefore $\Theta(\sigma)<\Theta(\eta)$. It is easy to see that the same argument also proves that $\Theta(\sigma)<\Theta(\eta)$ implies that $\sigma < \eta$ in the lexicographical order.

Consider now $\sigma<\eta$ and $x_\sigma$, $x_\eta$ defined as in eq. (\ref{eq-order10}). Define $k$ as before and suppose that $k<n$. Write $y=M_{\sigma_{k+1}} \circ \cdots \circ M_{\sigma_n} (0)$, $z = M_{\sigma_k}(y)$, so that $x_\sigma = M_{\sigma_1} \circ  \cdots \circ M_{\sigma_{k-1}} (z)$.
Observe that $y$ is less than, or equal to $M_1^{n-k}(0)=1-\frac{1}{n-k+1}$, so that $z \leq M_0(1-\frac{1}{n-k+1}) = \frac{n-k}{2n-2k+1} < \frac{1}{2}$.
Equivalently, write $u=M_{\eta_{k+1}} \circ \cdots \circ M_{\eta_n} (0)$, $v = M_{\eta_k}(u)$, so that $x_\eta = M_{\eta_1} \circ  \cdots \circ M_{\eta_{k-1}} (v)$. Now, $u \geq 0$, so that $v = M_{1}(u) \geq \frac{1}{2}$, and therefore $v > z$.
The map composition $M_{\eta_1} \circ  \cdots \circ M_{\eta_{k-1}}$ is the same as $M_{\sigma_1} \circ  \cdots \circ M_{\sigma_{k-1}}$, since $\sigma_j=\eta_j$ for $j<k$; being composed of strictly increasing maps is itself strictly increasing, so that $z<v$ implies $x_\sigma < x_\eta$.
It remains to consider the case $k=n$. In this case, $\sigma=\upsilon 0$, $\eta=\upsilon 1$, with $\upsilon \in \Sigma^{n-1}$. Therefore $x_\sigma=M_\upsilon(0)$, which is smaller than $x_\eta=M_\upsilon(\frac{1}{2})$.

Let us now prove the third statement of the lemma. When $n=0$ we have that $I_\emptyset = [0,1]$ and
$x_\emptyset = M_\emptyset(0) = 0$ (because $M_\emptyset$ is the identity); also $x_{\hat \sigma} = 1$, because $\emptyset$ is $1^0$, so that $x_{\hat \sigma} = 1$ by definition,  so that eq. (\ref{eq-order1b}) holds. When $n>0$, $I_\sigma= [M_\sigma(0),M_\sigma(1)] = [x_\sigma,M_\sigma(1)]$: we have to prove that
$M_\sigma(1) = x_{\hat \sigma}$. Clearly, when $\sigma=1^n$ $M_\sigma(1)=1$ and, by the definition above, $x_{\hat \sigma}=1$. Suppose that $M_\sigma(1) = x_{\hat \sigma}$ holds for any $\sigma \in \Sigma^n$. This is clearly true for $n=1$, since either $\sigma=0$, $\hat \sigma = 1$ and $M_0(1)=M_1(0)=x_1$, or $\sigma=1$, $M_1(1)=1$ and by definition $x_{\hat \sigma} =1$. Consider now $\sigma \in \Sigma^{n+1}$. Write $\sigma = \eta i$ with $\eta \in \Sigma^n$, $i=0,1$. In the first case
\[
 M_{\sigma}(1) = M_\eta M_0 (1) = M_\eta M_1 (0) = M_{\eta 1}(0) = x_{\eta 1}
 \]
and clearly $\eta 1 = \hat \sigma$. In the second case, suppose that $\eta \neq 1^n$, since the opposite instance means $\sigma=1^{n+1}$, which was treated above. Then, using the induction hypothesis and the fact that $M_0(0)=0$ we obtain
\[
 M_{\sigma}(1) = M_\eta M_1 (1) = M_\eta (1) = M_{\hat \eta}(0) =
 M_{\hat \eta} M_0 (0) = M_{\hat \eta 0}(0) = x_{\hat \eta 0}.
 \]
Since $\hat \sigma = \widehat{\eta 1} = \hat \eta 0$ the thesis follows. $\Box$

\begin{lemma}
Let $\Sigma^n$ be as in definition \ref{def-1} and $x_\sigma$, $I_\sigma$  be defined as in Lemma \ref{lem-2a}, eqs. (\ref{eq-order10}) and (\ref{eq-order1b}).
Then, for any $n \in \bfn$, $\sigma \in \Sigma^n$
\beq
?(x_\sigma) = \sum_{j=1}^{n} \sigma_j 2^{-j} = 2^n \Theta(\sigma)
\nuq{eq-value1}
and 
\beq \mu(I_\sigma) = 2^{-n}.
\nuq{eq-measi1}
\label{lem-2}
\end{lemma}
{\em Proof}.
Let us first prove eq. (\ref{eq-value1}).
From eq. (\ref{eq-mink1b}) it follows that $?(x_\sigma) = P_\sigma(0)$ for any $\sigma \in \Sigma$. Let us use induction again. For $n=0$ we have that $\sigma=\emptyset$ and eq. (\ref{eq-order1}) implies that $\Theta(\emptyset)=0=?(0)$. For $n=1$ we have that $x_0 = 0$ and $?(0)=0$; $x_1=\frac{1}{2}$ and $?(x_1)=\frac{1}{2}$, which again confirms eq. (\ref{eq-value1}). Next, suppose that eq. (\ref{eq-value1}) holds in $\Sigma^n$ and let us compute $?(x_{\sigma})$,  with $\sigma \in \Sigma^{n+1}$. Clearly, $\sigma = i \eta$, with $i=0$ or $i=1$, $\eta \in \Sigma^n$. Therefore,
\[
?(x_\sigma) =
?(x_{i\eta}) = P_i (?(x_\eta)) = P_i (\sum_{j=1}^{n} \eta_j 2^{-j}),
\]
Since $P_i(y)= i/2  + y/2$ we find
\[
?(x_{i\eta}) = i\; 2^{-1} + \sum_{j=1}^{n} \eta_j 2^{-j-1},
\]
which proves formula (\ref{eq-value1}).

Let us now compute $\mu(I_\sigma) = \mu([x_\sigma,x_{\hat \sigma}]) = ?(x_{\hat \sigma}) - ?(x_\sigma)$. When $n=0$, $\sigma=\emptyset$ we have that $I_\sigma=[0,1]$ so that $\mu(I_\sigma) = 1$.
When $\sigma \neq 1^n$ we can use eq. (\ref{eq-value1}), to obtain $?(x_{\hat \sigma}) - ?(x_\sigma) = 2^{-n}[\Theta(\hat \sigma)-\Theta(\sigma)] = 2^{-n}$.
If $\sigma = 1^n$ then $x_{\hat \sigma}=1$ and $?(x_{\hat \sigma}) - ?(x_\sigma) = ?(1) - ?(x_{1^n}) = 1 - 2^{-n}(2^n-1) = 2^{-n}$ and eq. (\ref{eq-measi1}) holds.$\Box$

\section{Stern--Brocot sequences and M\"obius IFS}
\label{sec-2}

In this section we demonstrate that the classical Stern--Brocot sequences \cite{stern,brocot,knuth} are realized as the boundary points of the M\"obius IFS partitions just described.

\begin{definition}
The Stern--Brocot sequence $B^n \subset \bfq$ is defined for any $n \in \bfn$ by induction: $B^0 = \{0,1\}$ and $B^{n+1}$ is the increasingly ordered union of $B^n$ and the set of mediants of consecutive terms of $B^n$. Recall that the mediant, or {\em Farey sum}, of two rational numbers written as irreducible fractions, is
\[
  \frac{p}{q} \oplus \frac{\hat p}{\hat q} = \frac{p + \hat p}{q + \hat q}.
\]
\end{definition}
Observe that the mediant of two numbers is intermediate between the two. Moreover, the definition implies that the cardinality of $B^n$ fulfills the rules $\#(B^0)=2$, $\#(B^{n+1}) = 2 \#(B^n) - 1$, whose solution is $\#(B^n)=2^n+1$. Therefore, the construction rule can be written as
\beq
B^n = \{x^n_0,x^n_1,x^n_2,\ldots,x^n_{2^n}\} \Rightarrow B^{n+1} = \{x^n_0,x^n_0\oplus  x^n_1,x^n_1,x^n_1\oplus  x^n_2,x^n_2,\ldots,x^n_{2^n}\}.
\nuq{eq-fa1}
The above equation also serves to introduce a symbolic notation for $B^n$. The next important lemma draws the relation between Stern--Brocot sequences and the partitions of $[0,1]$ generated by the M\"obius Iterated Function System (\ref{eq-mink1}).

\begin{lemma}
\label{lem-3}
Let $\Sigma^n$ be as in definition \ref{def-1} and $x_\sigma$, $I_\sigma$, for $\sigma \in \Sigma^n$, be defined as in Lemma \ref{lem-2a}, eqs. (\ref{eq-order10}) and (\ref{eq-order1b}). For any $n \in \bfn$ the increasingly ordered set  $\{\{x_{\sigma}, \sigma \in \Sigma^n\},1\}$ coincides with the $n$-th Stern-Brocot sequence $B^n$.
\end{lemma}

Observe that $\{\{x_{\sigma}, \sigma \in \Sigma^n\},1\}$ is the ordered set of extrema of the intervals $I_\sigma$, with $\sigma \in \Sigma^n$, which is increasingly ordered according to Lemma \ref{lem-2a}.
For $n=0$ one has $\{x_{\emptyset},1\}=\{0,1\}$, which can also be written as $B^0 = \{\frac{0}{1},\frac{1}{1}\}$. It is then enough to show that the induction property (\ref{eq-fa1}) holds for the set $\{\{x_{\sigma}, \sigma \in \Sigma^n\},1\}$.
Let $\sigma \in \Sigma^n$. Each $I_\sigma = [x_\sigma,x_{\hat \sigma}]$ splits into $I_{\sigma 0}$ and $I_{\sigma 1}$, as seen above in Lemma \ref{lem-1a}. Because of eq. (\ref{eq-cover}) the points $x_\sigma$ and $x_{\hat \sigma}$ of the $n$-th set also belong to the $n+1$-th set: in fact, they coincide with $x_{\sigma 0}$ and $x_{\hat \sigma 0}$. It remains to show that the intermediate point $x_{\sigma 1}$ is a rational that fulfills the Farey sum rule.
Remark that $M_\sigma$ is a M\"obius transformation in $PSL (2,\bfz)$, because both $M_i$, $i=0,1$ are such. Let
\beq
 M_\sigma(0) = x_\sigma = \frac{{p}}{{q}},
 \;
 M_\sigma(1) = x_{\hat{\sigma}} = \frac{\hat{p}}{\hat{q}},
 \nuq{eq-extre1}
where $p$ and $q$, $\hat p$ and $\hat q$ are relatively prime integers.
Observe that, when $n=0$,
\beq
 \Delta(\frac{{p}}{{q}},\frac{\hat{p}}{\hat{q}}) = \hat p {q} - \hat q {p} = 1.
 \nuq{eq-delta1}
Suppose that property (\ref{eq-delta1}) holds for $\sigma \in \Sigma^n$. Eqs. (\ref{eq-extre1}) and (\ref{eq-delta1}) imply that the M\"obius transformation $M_\sigma$ can be explicitly written as
\beq
M_\sigma : x \rightarrow \frac{(\hat p-{p})x +{p}}{(\hat q- {q})x +{q}}.
 \nuq{eq-delta2}
Therefore,
\[
x_{\sigma 1} = M_{\sigma 1}(0) = M_{\sigma} (\frac{1}{2}) =
\frac{p +\hat{p}}{q +\hat{q}},
\]
which proves that the Farey sum property holds. Using the last formula one also proves that eq. (\ref{eq-delta1}) {\em i.e.} $\Delta(\frac{p}{q},\frac{\hat{p}}{\hat{q}})=1$ holds between neighboring Stern--Brocot points of $B^{n+1}$. $\Box$

\section{Proof of regularity of Minkowski's question mark measure}
\label{sec-3}

After the preparatory work of the previous two sections, we now move to the core of the matter. In this section we first reformulate Criterion $\lambda^*$ in a way suitable for further developments. Next, we prove a key Lemma that links Criterion $\lambda^*$ to the IFS partitions of $[0,1]$ introduced before. We then state and prove a difficult result on the cardinality of ``large'' intervals of these partitions. Finally, we link all these results together in the proof of Theorem \ref{teo-1}.

Let us first observe that Criterion $\lambda^*$ is clearly
equivalent to the following requirement:
\begin{criterion}[$\lambda^*$]
\label{cri-one}
Define
\beq
\Lambda^{n}(\alpha) = \{ x \in [0,1] \mbox{ s.t. } \mu([x,x+\alpha/n]) \geq 2^{-n} \}.
\nuq{eq-reg00}
If for any $\alpha >0$, the Lebesgue measure of $\Lambda^{n}(\alpha)$ tends to one, when $n$ tends to infinity, then the positive Borel measure $\mu$ supported on $[0,1]$ is regular according to Ullman--Stahl--Totik.
\end{criterion}

Observe that, being Minkowski's question mark measure $\mu$ continuous, it is inessential whether the open or closed interval $[x,x+\alpha/n]$ is used in the above.
This rescaling of criterion $\lambda^*$ suits our purposes, for we can then rely on the structure of the M\"obius IFS / Stern--Brocot partitions of the interval $[0,1]$, which were described in Sections \ref{sec-1} and \ref{sec-2}. In fact, the previous theory will enter the argument through the fundamental lemma:
\begin{lemma}
Let $\alpha>0$, $n \in \bfn_+$ and let $B^n = \{x^n_j\}_{j=0,\ldots,2^n}$ be the $n$-th Stern--Brocot sequence. Define
\beq
A^n(\alpha) = \{ j \mbox{ s.t. } 0 \leq j \leq 2^n-2, \;  x^n_{j+2}-x^n_j \leq \alpha / n \}.
\nuq{eq-seta}
Then,
\beq
\bigcup_{j \in A^n(\alpha)} [x^n_j,x^n_{j+1}] \subset \Lambda^n(\alpha)
\nuq{eq-setb}
\label{lem-4}
\end{lemma}
{\em Proof.}
Because of Lemma \ref{lem-2} we have that
$
   ?(x^n_j) = j \; 2^{-n}.
$
Let $j \in A^n(\alpha)$, that is, the following inequality holds: $x^n_{j+2}-x^n_j \leq \alpha / n$. Then, for any $x \in [x^n_j,x^n_{j+1}]$ one has that $?(x) \leq ?(x^n_{j+1}) = (j+1) \; 2^{-n}$ and also, since $x+\alpha/n \geq x^n_{j+2}$, $?(x+\alpha/n) \geq ?(x^n_{j+2})=(j+2) \; 2^{-n}$. Therefore, for such $j$ and $x$, $\mu([x,x+\alpha/n]) = ?(x+\alpha/n) - ?(x) \geq 2^{-n}$, which means that $[x^n_j,x^n_{j+1}] \subset \Lambda^n(\alpha)$. $\Box$

We can now move on to a technical result, which needs care in notation.
Recall that the intervals $[x^n_j,x^n_{j+1}]$ in the above lemma correspond to $I_\sigma$, when $\sigma \in \Sigma^n$ is such that $\Theta(\sigma)=j$. We shall repeatedly pass from the integer order to the symbolic representation and back: unless otherwise stated, we always assume that $|\sigma|=n$, $\Theta(\sigma)=j$ and we write $I_\sigma = [x^n_{j}, x^n_{j+1}] = [ \frac{p^n_j}{q^n_j},\frac{p^n_{j+1}}{q^n_{j+1}}]$, with $p^n_j$ and $q^n_j$, $p^n_{j+1}$ and $q^n_{j+1}$ relatively prime integers. For the same object, we also use the notation
$I_\sigma = [x_\sigma,x_{\hat \sigma}] = [\frac{p}{q},\frac{\hat p}{\hat q}]$.
We let $\lambda$ denote Lebesgue measure, $\lambda(I_\sigma) = x_{\hat \sigma} - x_\sigma$ and finally, for simplicity of notation, we also set $a=1/\alpha$.

The coming proposition deals with ``large'' intervals $L^n(\alpha)$ ($L$ for large) in the IFS partitions:

\begin{proposition}
\label{prop-1}
Define the set $L^n(\alpha) \subset \Sigma^n$, for $n \in \bfn$, as
\beq
 L^n(\alpha) = \{ \sigma \in \Sigma^n \mbox{ s.t. }
  \lambda(I_\sigma) =  x^n_{j+1} - x^n_j \geq \frac{\alpha}{n}, \; \mbox{ where } j = \Theta(\sigma) \}.
  \nuq{eq-cna}
Then, for any $\alpha>0$ cardinality of $L^n(\alpha)$ is uniformly superiorly bounded: there exists $l(\alpha) \in \bfn$ so that
\beq
 \# (L^n(\alpha)) \leq l(\alpha), \; \forall n \in \bfn
\nuq{eq-cna2}
\end{proposition}
{\em Proof.} To prove this proposition we will proceed through several steps, some of which contain results that can be considered as sublemmas in their own right. We start by providing some definitions used in the proof.

{\em Definitions of useful sets.}

For any value of $\alpha>0$ we single out a class of rational in $[0,1]$, $\bfq_\alpha$, by considering all irreducible fractions with denominator smaller than $1/\sqrt{\alpha}$:
\beq
 \bfq_\alpha = \{ x \in \bfq  \cap [0,1] \mbox{ s.t. } x = \frac{p}{q}, \; p,q \in \bfn, \; p \perp q \mbox{ and } q^2 < \frac{1}{\alpha}\}.
 \nuq{eq-qu1}
We will show that the set $\bfq_\alpha$ determines $L^n(\alpha)$: the finite cardinality of the former will serve to bound the cardinality of the latter, in a way independent of $n$.

We also define the set of ``small intervals/words'': it is the complementary set of $L^n(\alpha)$ in $\Sigma^n$, which we call $S^n(\alpha)$ ($S$ for {\em small}):
\beq
  S^n(\alpha) = \{ \sigma \in \Sigma^n \mbox{ s.t. }
  \lambda(I_\sigma) < \alpha/n \}.
  \nuq{eq-dna}
The property of Farey fractions, eq. (\ref{eq-delta1}), imply that
\beq
  \lambda(I_\sigma) = x^n_{j+1} - x^n_j  = \frac{1}{q^n_{j+1} q^n_j}.
  \nuq{eq-cna6}
This permits to assess the useful condition
\beq
 \sigma \in S^n(\alpha) \Longleftrightarrow q^n_{j+1} q^n_j > a n, \mbox{ where } n=|\sigma|,  \; j=\Theta(\sigma).
 \nuq{eq-new}

Finally, we define a subset of $\Sigma$ by requiring that neither extremum of $I_\sigma$ belongs to $\bfq_\alpha$:
\beq
 {\cal E} = \{ \sigma \in \Sigma \mbox{ s.t. } q^n_j  > \sqrt{a}, \; q^n_{j+1} > \sqrt{a} \mbox{ when } n=|\sigma|, j = \Theta(\sigma) \}.
 \nuq{eq-es1}
Notice that in the above we do {\em not} require ${\cal E}$ to be a subset of $\Sigma^n$, but rather of the full set $\Sigma$: this is necessary to study different IFS partitions.

{\em Sublemmas}

Following the above notations and definitions we can prove a series of implications, or sublemmas.
Let us first show that
\beq
\sigma \in {\cal E} \Rightarrow  \sigma \eta \in {\cal E}, \; \forall \eta \in \Sigma.
\nuq{eq-cale1}

In fact, denominators of the endpoints of $I_{\sigma \eta}$ belong to the Stern--Brocot sequence $B^{|\sigma|+|\eta|}$. Let $I_{\sigma}= [\frac{p}{q},\frac{\hat p}{\hat q}]$. When $|\eta|=1$, because of the construction rule (\ref{eq-fa1}), the denominators of the endpoints of $I_{\sigma \eta}$ are $\{q,q+\hat q,\hat q\}$, which are all larger than $\sqrt{a}$. Induction extends the result to general $\eta \in \Sigma$. In words, this means that the class ${\cal E}$ is stable under successive partitions. $\triangle$

A second implication considers the case of words which are in ${\cal E}$ and at the same time are associated with ``small'' intervals. This class is also stable under successive partitions:
\beq
\sigma \in [{\cal E} \cap S^{|\sigma|}(\alpha)] \Rightarrow \sigma \eta \in [{\cal E} \cap S^{|\sigma|+|\eta|}(\alpha)], \; \forall \eta \in \Sigma.
\nuq{eq-eta1}
The implication regarding $\cal E$ has just been proven, relation (\ref{eq-cale1}). In addition, letting $I_\sigma = [\frac{p}{q},\frac{\hat p}{\hat q}]$, the l.h.s. of (\ref{eq-eta1}) means that $q^2 >a$, $\hat{q}^2 >a$ and $q \hat{q} >a |\sigma|$. Then, $I_{\sigma 0} = [\frac{p}{q},\frac{\hat p+p}{\hat q+q}]$ and
\[
  q ({\hat q+q}) = q \hat q + q^2 > a |\sigma|  +  a =
   a (|\sigma| + 1) = a |\sigma 0| .
\]
A similar estimate clearly holds for $I_{\sigma 1}$: by induction, this proves (\ref{eq-eta1}). $\triangle$

Suppose now that $\sigma \in {\cal E}$, but we do not require that  $\lambda(I_\sigma)< \alpha/|\sigma|$, {\em i.e. } $\sigma$ may {\em not} belong to $S^{|\sigma|}(\alpha)$, a case which often occurs. Then, we can prove that there is a subdivision of $I_\sigma$  whose intervals are all smaller than the threshold in (\ref{eq-dna}), {\em i.e.} the words $\sigma \eta$ in this subdivision belong to $S^{|\sigma|+|\eta|}(\alpha)$:
\beq
\sigma \in {\cal E} \Rightarrow \exists \; k_1(\sigma) \in \bfn \mbox{ s.t. } \sigma \eta \in [{\cal E} \cap S^{|\sigma|+|\eta|}(\alpha)],  \;
\forall \eta \in \Sigma, \;|\eta| \geq k_1(\sigma).
\nuq{eq-eta3}
The proof of this implication is rather long.
Again, the part regarding $\cal E$ has been proven above, eq. (\ref{eq-cale1}). Let us use again the notation
$I_\sigma = [\frac{p}{q},\frac{\hat p}{\hat q}]$.
Suppose that $q < {\hat q}$ without loss of generality: the opposite case can be dealt with similarly, by replacing the symbols $q$ with $\hat q$ and $0$ with $1$ in the following.
In the former case, among all intervals $I_{\sigma \eta}$, with $|\eta|=k$, the largest is $I_{\sigma 0^k}$, as we are going to prove in the next two paragraphs.

Since $I_{\sigma \eta} = M_\sigma(I_\eta)$, let us first study the intervals $I_\eta$ and prove that $\lambda(I_{\eta}) \leq \lambda(I_{0^{k}})$ for any $\eta \in \Sigma^k$ and for any $k \in \bfn$. This is clearly true for $k=0$ and $k=1$ by direct inspection. Suppose that it holds true for a certain $k$ and let us consider symbolic words of length $k+1$. Since the associated intervals are symmetric around $1/2$, it is sufficient to study the set $\{I_{0 \eta}$, $\eta \in \Sigma^k\}$. We now want to use the mean value theorem: the equality $I_{0 \eta} = M_0(I_{\eta})$ holds by definition, so that $\lambda(I_{0 \eta}) = M'_0(z_{\eta}) \lambda(I_{\eta})$, where $z_{\eta}$ is a point in $I_{\eta}$. Now, $M'(x)=1/(x+1)^2$, so that $M'(x)$ is strictly decreasing on $[0,1]$. Then,
since for any $\eta \in \Sigma^k$ we have that $I_\eta \geq I_{0^k}$, by the ordering of intervals in Lemma \ref{lem-2a}, and because $\lambda(I_\eta) \leq \lambda(I_{0^k})$ by the induction hypothesis, $\lambda(I_{0 \eta}) = M'_0(z_{\eta}) \lambda(I_{\eta}) \leq M'_0(z_{0^k}) \lambda(I_{0^k}) = \lambda(I_{0^{k+1}})$, which completes the induction proof. We can also explicitly compute
\beq
\lambda(I_{0^k}) = \lambda(I_{1^k}) = \frac{1}{k+1},
\nuq{eq-large}
which will be useful below.

Recall now that
$I_{\sigma \eta} = M_\sigma(I_\eta)$.
When $I_\sigma = [\frac{p}{q},\frac{\hat p}{\hat q}]$ the M\"obius transformation $M_\sigma$ is given by eq. (\ref{eq-delta2}), so that $M'_\sigma(x) = [(\hat q - q)x+q]^{-2}$, using also eq. (\ref{eq-delta1}). Since $q < \hat q$, $M'_\sigma(x)$ is decreasing on $[0,1]$. Then, using again the mean value theorem, $\lambda(I_{\sigma \eta}) = M'_\sigma(z_\eta) \lambda(I_\eta)$, with $z_\eta \in I_\eta$, we conclude that $\lambda(I_{\sigma \eta}) \leq \lambda(I_{\sigma 0^k})$ for any $\eta \in \Sigma^k$.

The length of $I_{\sigma 0^k}$ can be easily computed from the explicit representation
\[
 I_{\sigma 0^k} = [\frac{p}{q},\frac{\hat p +kp}{\hat q + kq}],
\]
which yields
\[
\lambda(I_{\sigma 0^k})^{-1} = q (\hat q + kq) = q \hat q + k q^2.
\]
Call $k_1(\sigma)$ the least $k$ such that $k > (na- q \hat q)/(q^2-a)$, where $n=|\sigma|$. Since $q^2-a>0$ this implies that, for $k \geq k_1(\sigma)$,
\beq
 q \hat q + k q^2 > a (n+k).
 \nuq{eq-eta4}
Inequality (\ref{eq-eta4}) and the above reasoning imply that
\[
  \lambda(I_{\sigma \eta}) \leq \lambda (I_{\sigma 0^k})\leq \frac{\alpha}{n+k}
  \]
for all $\eta \in \Sigma^k$, $k \geq k_1(\sigma)$, which proves (\ref{eq-eta3}). $\triangle$

We have just proven that that starting from a word in ${\cal E}$ and taking successive partitions we always end up in ${\cal E} \cap S^{m}(\alpha)$ for all $m \in \bfn$ larger than a certain value. We now need to examine the fate of intervals which do {\em not} necessarily belong to ${\cal E}$.

Let us therefore take a general $\sigma \in \Sigma$, not necessarily in ${\cal E}$ and consider the associated interval $I_\sigma = [\frac{p}{q},\frac{\hat p}{\hat q}]$, where either $q^2 \leq a$ or ${\hat q}^2 \leq a$ may happen.
Nonetheless, we can prove that
\beq
\forall \sigma \in \Sigma, \; \; \exists \; k_2(\sigma) \in \bfn \mbox{ s.t. } \sigma 0^k 1 \in [{\cal E} \cap S^{|\sigma|+k+1}(\alpha)], \; \forall k \geq k_2(\sigma).
\nuq{eq-eta6}
By direct computation one gets
\[
 I_{\sigma 0^k 1} = [\frac{\hat p +(k+1)p}{\hat q + (k+1)q},\frac{\hat p +kp}{\hat q + kq}].
\]
Observe that these intervals approach $x_\sigma=\frac{p}{q}$ when $k$ grows. Actually, $I_{\sigma 0^k 1}$ is the second interval to the right of $\frac{p}{q}$ in the family $\{I_\eta, \; |\eta| = |\sigma| + k + 1\}$.
It is clear that for sufficiently large $k$ the squares of both denominators are larger than $a$, so that $\sigma 0^k 1 \in {\cal E}$. Moreover, since
\[
\lambda(I_{\sigma 0^k 1})^{-1} =   [\hat q + (k+1)q] (\hat q + kq)= k^2 q^2 + k q^2 + (2k+1)\hat q  q + \hat q ^2 .
\]
it is also clear that, for sufficiently large $k$, the r.h.s. of the above is larger than $a(n+k+1)$, so that $\sigma 0^k 1 \in S^{n+k+1}$, where $n=|\sigma|$. This proves (\ref{eq-eta6}). $\triangle$

While fully general, we will use the above property for $\sigma$ such that $q^2 \leq a$.
We also need a symmetrical property, to be used when $\hat q^2 \leq a$. Using the same technique we easily show that
\beq
\forall \sigma \in \Sigma, \; \; \exists \; k_3({\sigma}) \in \bfn \mbox{ s.t. } \sigma 1^k 0 \in [{\cal E} \cap S^{|\sigma|+k+1}(\alpha)], \; \forall k \geq k_3( \sigma).
\nuq{eq-eta6b}
Here $I_{\sigma 1^k 0}$ is the second interval to the left of $\frac{\hat p}{\hat q}$ in the family $\{I_\eta, \; |\eta| = |\sigma| + k + 1\}$. $\triangle$

We now go through a series of three levels $n=n_1,n_2,n_3$ at which we study the partitions $\Sigma^n$. \\
{\em First level, $n_1$, when all elements of $\bfq_\alpha$ appear in $B^n$.}\\ %
Consider the set $\bfq_\alpha$, whose cardinality is obviously finite. 
For $\zeta \in \bfq_\alpha$ let $n(\zeta)$ be the least $n$ such that $\zeta \in B^n$. This number is sometimes called the {\em depth} of $\zeta$ in the Stern--Brocot tree. Note that it exists for any $\zeta \in \bfq \cap [0,1]$ (see {\em e.g.} \cite{knuth}), although we have not proven it here. Nonetheless, our proof does {\em not} require this result, for we only need to consider values $\zeta \in \bfq_\alpha$ that appear in Stern--Brocot sequences. 
Since $\zeta \in B^n$ implies that $\zeta \in B^m$ for any $m \geq n$, there exists a value $\bar n = \max \{ n(\zeta), \zeta \in \bfq_\alpha \}$ such that $\zeta \in B^{n}$ for all $\zeta \in \bfq_\alpha$ and for all $n \geq \bar n$. Let now $n_1  = \bar n+ 1$, so that if $x^{n_1}_j \in \bfq_\alpha$ then $j$ is even: in fact, for any $\zeta \in \bfq_\alpha$ there exists $j \in \{0,\ldots,2^{\bar n}\}$ such that $\zeta=x^{\bar n}_j$. At the next level, we have $\zeta=x^{\bar n+1}_{2j}$, so that values of $\bfq_\alpha$ appear with even indices in $B^{n_1}$. Since this exhausts all elements of $\bfq_\alpha$, no odd index terms of $B^{n_1}$ belong to this set.

Consider now the set $F$ of words in  $\Sigma^{n_1}$, such that one endpoint of $I_\sigma$ belongs to $\bfq_\alpha$. Because of what we have just proven, no more than one endpoint of a single interval can belong to $\bfq_\alpha$. Part these words in two groups, according to whether the left or the right endpoint of $I_\sigma$ lie in $\bfq_\alpha$:
\[
  F_l = \{ \sigma \in \Sigma^{n_1} \mbox{ s.t. } x_\sigma \in \bfq_\alpha \},
\]
\[
  F_r = \{ \sigma \in \Sigma^{n_1} \mbox{ s.t. } x_{\hat \sigma} \in \bfq_\alpha \}.
\]
Apply now sublemmas (\ref{eq-eta6}) and
(\ref{eq-eta6b}) to define $K_2(F) = \max \{k_2(\sigma), \sigma \in F_l \}$ and
$K_3(F) = \max \{k_3(\sigma), \sigma \in F_r \}$. Let $\kappa = \max \{K_2(F),K_3(F)\} + 1$. This defines the second level, $n_2$.

{\em Second level, $n_2=n_1+\kappa$, when properties (\ref{eq-eta6}) and
(\ref{eq-eta6b}) are realized for all words $\sigma$ in $F \subset \Sigma^{n_1}$}

Among all words of length $n_2= n_1 + \kappa$ we start by considering those that originate from a word $\sigma \in F_l$. They are written as $\sigma \eta$, where $\eta$ is any word in $\Sigma^{\kappa}$. All of these belong to ${\cal E}$, except for $\sigma 0^{\kappa}$.  Equally, when $\sigma \in F_r$, the words $\sigma \eta$, where $\eta$ is any word in $\Sigma^{\kappa}$, belong to ${\cal E}$, except for $\sigma 1^{\kappa}$. Because of the argument above, these two cases yield all words of $\Sigma^{n_2}$ that are {\em not} in $\cal E$. We can therefore write $\Sigma^{n_2}$ as the union of three sets that are pairwise disjoint:
\beq
  \Sigma^{n_2} = ({\cal E} \cap \Sigma^{n_2}) \uplus \{\sigma 0^{\kappa}, \; \sigma \in F_l\}
\uplus \{\sigma 1^{\kappa}, \; \sigma \in F_r\}.
\nuq{eq-disjo}
Consider now the first set in the disjoint union above, call it $E  = {\cal E} \cap \Sigma^{n_2}$. It corresponds to intervals $I_\sigma$, with $\sigma \in \Sigma^{n_2}$, such that neither endpoint of $I_\sigma$ belongs to $\bfq_\alpha$. Two cases are possible: small and large intervals.
\[
  E_s = E \cap S^{n_2}(\alpha), \; E_l = E \cap L^{n_2}(\alpha).
\]
In the first case, that is, $\sigma \in E_s$, $I_{\sigma \eta}$ belongs to ${\cal E} \cap S^{n_2 + k}(\alpha)$ for any $k \geq 0$, $\eta \in \Sigma^k$, in force of (\ref{eq-eta1}). In the second case, $\sigma \in E_l$, we use (\ref{eq-eta3}): for any such $\sigma \in E_l$  there exists $k_1(\sigma)$ such that $\sigma \eta \in {\cal E} \cap S^{n_2 +k}(\alpha)$ for any $k \geq k_1(\sigma)$, $\eta \in \Sigma^k$.
Since the cardinality of $E_l$ is finite, there exists a maximum of $\{k_1(\sigma), \; \sigma \in E_l \}$: call it $K_1(E)$. This defines a new level, $n_3=n_2+K_1(E)$.

{\em Third level, $n_3$, where we make the final separation between small and large intervals.}


We have just proven that for any $n \geq n_3$ the words $\sigma \eta$, with $\sigma \in {\cal E} \cap \Sigma^{n_2}$ and $\eta \in \Sigma^{n-n_3}$ belong to ${\cal E} \cap S^{n}(\alpha)$.
It remains to consider words in $\Sigma^{n}$, with $n \geq n_3$, which originate from
$\{\sigma 0^{\kappa}, \; \sigma \in F_l\}$ and $\{\sigma 1^{\kappa}, \; \sigma \in F_r\}$: recall that these sets are included in $\Sigma^{n_2}$; we need to consider their offsprings. Let us show how to proceed by induction. Consider the first case and start from the word $\sigma 0^{\kappa}$, with $\sigma \in F_l$. Its partition yields the two words $\sigma 0^{\kappa+1}$ and $\sigma 0^{\kappa} 1$. Because of (\ref{eq-eta6}) and because $\kappa > k_2(\sigma)$, the latter belongs both to $\cal E$ and to $S^{n_2+1}(\alpha)$: it is associated with a small interval. We can now apply (\ref{eq-eta1}) to show that all of its offsprings $\sigma 0^{\kappa} 1 \eta$, with $\eta \in \Sigma^{m}$, for every $m \in \bfn$ belong to ${\cal E} \cap S^{n_2+1+m}(\alpha)$. We iterate the procedure on $\sigma 0^{\kappa+1}$ and so on, which proves that for any $m \in \bfn$ all words of the kind
\[
  \{\sigma 0^{\kappa} \eta, \; \sigma \in F_l, \eta \in \Sigma^{m} \}
\]
belong to $S^{n_2+m}(\alpha)$, except possibly for $\{\sigma 0^{\kappa+m},  \; \sigma \in F_l\}$. Similarly, we prove that all words of the kind
\[
  \{\sigma 1^{\kappa} \eta, \; \sigma \in F_r, \eta \in \Sigma^{m} \}
\]
belong to $S^{n_2+m}(\alpha)$, except possibly for $\{\sigma 1^{\kappa+m},  \; \sigma \in F_r\}$.

{\em Conclusion}

The above classification of intervals shows that for all $n \geq n_3$
\beq
  \Sigma^n = S^n(\alpha) \cup \{\sigma 0^{n-n_1}, \; \sigma \in F_l \} \cup \{\sigma 1^{n-n_1}, \; \sigma \in F_r \},
\nuq{eq-fina1}
so that
\beq
  L^n(\alpha) \subset   \{\sigma 0^{n-n_1}, \; \sigma \in F_l \} \cup \{\sigma 1^{n-n_1}, \; \sigma \in F_r \}.
\nuq{eq-fina2}
The cardinality of $L^n(\alpha)$ in eq. (\ref{eq-fina2}) is less than the sum of the cardinalities of the two sets at r.h.s., which is obviously finite and independent of $n$, for all $n \geq n_3$. Clearly the maximum cardinality of $L^n(\alpha)$ for $n < n_3$ is also finite, so that we obtain the thesis. $\Box$
\begin{remark}
One can also prove that, for large $n$, all words in $\{\sigma 0^{n-n_1}, \; \sigma \in F_l \}$ and $\{\sigma 1^{n-n_1}, \; \sigma \in F_r \}$ do belong to $L^n(\alpha)$. From this one could therefore estimate precisely the Lebesgue measure of the set $\Lambda^n(\alpha)$ in Criterion \ref{cri-one}.
\end{remark}

We can now finally prove the main theorem of this paper: \\
{\bf Theorem 1} {\em Minkowski's question mark measure satisfies criterion $\lambda^*$ and hence is regular.}
%
{\em Proof.}
Let $\lambda$ denote Lebesgue measure. The theorem will be proven if we show that, for any $\alpha>0$, the Lebesgue measure of the complement of $\Lambda^n(\alpha)$, call it $\ovl{\Lambda}^n(\alpha)$, tends to zero, when $n$ tends to infinity.
Also as a notation, let $\ovl{A}^n(\alpha)$ the complementary set of $A^n(\alpha)$ in $\{0,\ldots,2^n-1\}$.
Because of Lemma \ref{lem-1a} and Lemma \ref{lem-4}, eq. (\ref{eq-setb}), we have that
\[
 \ovl{\Lambda}^n(\alpha) \subset \ovl{\bigcup_{j \in A^n(\alpha)} [x^n_j,x^n_{j+1}]} =  \bigcup_{j \in \ovl{A}^n(\alpha)} [x^n_j,x^n_{j+1}],
\]
where the last equality holds modulo a set of zero Lebesgue measure (which consists of interval endpoints). Therefore
\beq
\lambda(\ovl{\Lambda}^n(\alpha))\leq
\lambda(\bigcup_{j \in \ovl{A}^n(\alpha)} [x^n_j,x^n_{j+1}]) =
\sum _{j \in \ovl{A}^n(\alpha)}\lambda([x^n_j,x^n_{j+1}]) =
\sum _{j \in \ovl{A}^n(\alpha)} x^n_{j+1}- x^n_j.
\nuq{eq-prova01}
Observe now that
$ x^n_{j+1} - x^n_j \leq \frac{1}{n+1}$
as it was proven in Proposition \ref{prop-1}, eq. (\ref{eq-large}). We can therefore write
\beq
\lambda(\ovl{\Lambda}^n(\alpha))\leq \frac{1}{n+1}
\# ( \ovl{A}^n(\alpha)),
\nuq{eq-prova02}
where $\#(\ovl{A}^n(\alpha))$ indicates the cardinality of the set  $\ovl{A}^n(\alpha)$ that we need now to estimate.

First, because $\ovl{A}^n(\alpha)$ is the complementary set of $A^n(\alpha)$ in $\{0,\ldots,2^n-1\}$, it can be explicitly written as follows
\[
 \ovl{A}^n(\alpha) = \{ j \in \{0,\ldots,2^{n}-2 \}  \mbox{ s.t. }
  x^n_{j+2} - x^n_j \geq \frac{\alpha}{n} \} \cup \{2^n-1 \}.
  \]
Next, observe that  $x^n_{j+2} - x^n_j \geq \frac{\alpha}{n}$ implies that either $x^n_{j+1} - x^n_j \geq \frac{\alpha}{2n}$ or $x^n_{j+2} - x^n_{j+1} \geq \frac{\alpha}{2n}$, so that
\[
 \ovl{A}^n(\alpha) \subset \bigcup_{k=0,1} \{ j \in \{0,\ldots,2^{n}-2 \} \mbox{ s.t. }
  x^n_{j+1+k} - x^n_{j+k} \geq \frac{\alpha}{2n} \}
   \cup \{2^n-1 \}.
  \]
The cardinality of $\ovl{A}^n(\alpha)$ is less than, or equal to, the sum of the cardinality of the three sets at r.h.s. of the above equation. If we now recall the definition of the set $L^n(\alpha)$, eq. (\ref{eq-cna}), we see that the cardinality of the first two sets in the r.h.s. above is less than, or equal to, the cardinality of $L^n(\alpha/2)$, so that
\beq
\#(\ovl{A}^n(\alpha)) \leq 2 \#(L^n(\frac{\alpha}{2})) + 1.
\nuq{eq-prova3}
Using now the fundamental result of Proposition \ref{prop-1}, we obtain
 \beq
 \lambda(\ovl{\Lambda}^n(\alpha))\leq \frac{2}{n+1} l(\frac{\alpha}{2}) + \frac{1}{n+1},
\nuq{eq-prova4}
which tends to zero as $n$ tends to infinity and therefore it proves the thesis. $\Box$

{\bf Acknowledgements} \\
Studying abstract mathematical models is fun especially when it provides matter of entertainment with friends: I acknowledge fruitful and enjoyable discussions on the problems discussed herein with Walter Van Assche and Roberto Artuso.


\end{document}